\documentclass[11pt]{article}

\usepackage{amsmath, amssymb, amscd}
\usepackage[matrix,arrow]{xy}

\def\eqref#1{(\ref{#1})}
\newcommand{\goth}{\frak}
\newcommand{\g}{{\frak g}}
\newcommand{\gotht}{{\frak t}}
\newcommand{\arrow}{{\:\longrightarrow\:}}
\newcommand{\Z}{{\Bbb Z}}
\newcommand{\C}{{\Bbb C}}
\newcommand{\R}{{\Bbb R}}
\newcommand{\Q}{{\Bbb Q}}

\def\1{\sqrt{-1}\:}
\newcommand{\restrict}[1]{{\left|_{{\phantom{|}\!\!}_{#1}}\right.}}
\newcommand{\cntrct}                % contraction with a vector field
{\hspace{2pt}\raisebox{1pt}{\text{$\lrcorner$}}\hspace{2pt}}

\newcommand{\calo}{{\cal O}}

\renewcommand{\phi}{\varphi}
\renewcommand{\epsilon}{\varepsilon}
\renewcommand{\geq}{\geqslant}

\newcommand{\Tot}{\operatorname{Tot}}

\newcommand{\Aut}{\operatorname{Aut}}

\newcommand{\Lie}{\operatorname{Lie}}

\newcommand{\comment}[1]{{}}

\def\blacksquare{\hbox{\vrule width 4pt height 4pt depth 0pt}}
\def\endproof{\blacksquare}

\makeatletter

\@ifundefined{Bbb}
     {\newcommand{\Bbb}[1]{{\mathbb #1}}}%
{}%     {\edef\Bbb#1{{\Bbb #1}}}

\newcommand{\ps@verbit}{%
  \renewcommand{\@oddhead}{%
          \scriptsize
          {Positive toric fibrations}
          \hfil\tiny {M. Verbitsky, March 6, 2007 }}
  \renewcommand{\@evenhead}{\@oddhead}
  \renewcommand{\@oddfoot}{\hfil\thepage\hfil}
  \renewcommand{\@evenfoot}{\@oddfoot}}
 
\newcounter{Mycounter}[section]
\newcounter{lemma}[section]
\renewcommand{\thelemma}{\noindent{Lemma \thesection.\arabic{lemma}}}
\newcommand{\lemma}{%
     \setcounter{lemma}{\value{Mycounter}}
     \refstepcounter{lemma}
     \stepcounter{Mycounter}
     {\bf \thelemma:\ }}

\newcounter{claim}[section]
\renewcommand{\theclaim}{\noindent{Claim \thesection.\arabic{claim}}}
\newcommand{\claim}{%
     \setcounter{claim}{\value{Mycounter}}
     \refstepcounter{claim}
     \stepcounter{Mycounter}
     {\bf \theclaim:\ }}

\newcounter{sublemma}[section]
\newcounter{corollary}[section]
\newcounter{theorem}[section]
\renewcommand{\thetheorem}{\noindent{Theorem \thesection.\arabic{theorem}}}
\newcommand{\theorem}{%
     \setcounter{theorem}{\value{Mycounter}}
     \refstepcounter{theorem}
     \stepcounter{Mycounter}
     {\bf \thetheorem:\ }}

\newcounter{conjecture}[section]
\newcounter{proposition}[section]
\renewcommand{\theproposition}
       {\noindent{Proposition \thesection.\arabic{proposition}}}
\newcommand{\proposition}{%
     \setcounter{proposition}{\value{Mycounter}}
     \refstepcounter{proposition}
     \stepcounter{Mycounter}
     {\bf \theproposition:\ }}

\newcounter{definition}[section]
\renewcommand{\thedefinition}
       {\noindent{Definition~\thesection.\arabic{definition}}}
\newcommand{\definition}{%
     \setcounter{definition}{\value{Mycounter}}
     \refstepcounter{definition}
     \stepcounter{Mycounter}
     {\bf \thedefinition:\ }}

\newcounter{example}[section]
\newcounter{remark}[section]
\renewcommand{\theremark}{\noindent{Remark \thesection.\arabic{remark}}}
\newcommand{\remark}{%
     \setcounter{remark}{\value{Mycounter}}
     \refstepcounter{remark}
     \stepcounter{Mycounter}
     {\bf \theremark:\ }}

\newcounter{problem}[section]
\newcounter{question}[section]
\@addtoreset{equation}{section}
\@addtoreset{footnote}{section}
\makeatother

\begin{document}

%%%%%%%%%%%%%%%%%%%%%%%%%%%%%%%%%%%%%%%%%%%%%%%%%%%%%%%%%%%%
\begin{center}
{\LARGE\bf
Positive toric fibrations
}
%%%%%%%%%%%%%%%%%%%%%%%%%%%%%%%%%%%%%%%%%%%%%%%%%%%%%%%%%%%%
\\[4mm]
Misha Verbitsky\footnote{Misha Verbitsky is 
an EPSRC advanced fellow 
supported by  EPSRC grant  GR/R77773/01}
\\[4mm]

{\tt verbit@maths.gla.ac.uk, \ \  verbit@mccme.ru}
\end{center}

%%%%%%%%%%%%%%%%%%%%%%%%%%%%%%%%%%%%%%%%%%%%%%%%
{\small 
\hspace{0.15\linewidth}
\begin{minipage}[t]{0.7\linewidth}
{\bf Abstract} \\
A principal toric bundle $M$ is a complex manifold
equipped with a free holomorphic action of a compact
complex torus $T$. Such a manifold is fibered over
$M/T$, with fiber $T$. We discuss the notion of positivity
in fiber bundles and define positive toric bundles.
Given an irreducible complex subvariety $X\subset M$
of a positive principal toric bundle, we show that
either $X$ is $T$-invariant, or it lies in an orbit
of $T$-action. For principal elliptic bundles,
this theorem is known (\cite{_Verbitsky:Sta_Elli_}).
As follows from Borel-Remmert-Tits theorem,
any simply connected compact homogeneous 
complex manifold  is a principal toric bundle. 
We show that compact Lie groups with left-invariant 
complex structure $I$ are positive toric bundles,
if $I$ is generic. Other examples of positive
toric bundles are discussed.
\end{minipage}
}
%%%%%%%%%%%%%%%%%%%%%%%%%%%%%%%%%%%%%%%%%%%%%%%%

{
\small
\tableofcontents
}
%%%%%%%%%%%%%%%%%%%%%%%%%%%%%%%%%%%%%%%%%%%%%%%%

%%%%%%%%%%%%%%%%%%%%%%%%%%%%%%%%%%%%%%%%%%%%%%%%%%%%%%%%%%%%

\section*{Introduction}
\label{_Intro_Section_}

%%%%%%%%%%%%%%%%%%%%%%%%%%%%%%%%%%%%%%%%%%%%%%%%%%%%%%%%%%%%

Principal toric bundles were always a source of examples
in complex geometry (classical Hopf manifolds and Calabi-Eckmann
spaces come to mind), but the recent progress in string theory
in interaction with algebraic geometry puts a new light on
this classical subject. Since a discovery of the hypercomplex
structures on compact Lie groups and other homogeneous spaces
by the physicists Spindel et al (\cite{_SSTvP_}) and then Joyce 
(\cite{_Joyce_}), the complex geometry of homogeneous spaces
became a focus of much research, being a target space
of supersymmetric string theories with rich
supersymmetry. However, the algebraic geometry
of these spaces is not much understood. This paper
is an attempt to fill this gap, bringing a number of
results about subvarieties of principal
toric bundles, in particular - compact 
homogeneous compex manifolds with abelian 
fundamental group.

A principal toric bundle $M$ is a complex manifold
equipped with a free holomorphic action of a compact
complex torus $T$. Such a manifold is fibered over
$M/T$, with fiber $T$. In the sequel, we discuss 
 positivity in fiber bundles and define positive toric bundles,
generalizing a notion developed for
 elliptic bundles in \cite{_Verbitsky:Sta_Elli_}.

Many examples of principal toric bundles appear
from a theory of complex homogeneous manifolds.
It is known since 1960-es that a simply connected
compact complex homogeneous manifold is a principal
toric fibration over a rational projective homogeneous
manifold (this is known as {\bf the Tits fibration}).
 In fact, this result holds for a complex homogeneous manifold
with abelian fundamental group 
(\ref{_Borel_Remmert_Tits_genera_Theorem_}).

A principal toric fibration obtained this
way from a general complex structure on a Lie group
is always positive (\ref{_Samelson_positive_Theorem_}). 

We study complex subvarieties of positive principal 
toric fibrations.
Given an irreducible complex subvariety $X\subset M$
of a positive principal toric bundle, we show that
either $X$ is $T$-invariant, or it lies in an orbit
of $T$-action (\ref{_Q_pos_main_Theorem_}).

This theorem is elementary for $\dim_\C T=1$,
but for $\dim_\C T>1$ the proof uses highly non-trivial
results on deformation theory and the
spaces of Fujiki class C, which are due to Fujiki and 
Demailly-Paun.

%%%%%%%%%%%%%%%%%%%%%%%%%%%%%%%%%%%%%%%%%%%%%%%%%%%%%%%%%%%%

\section{Positivity in fiber bundles}
\label{_Intro_Section_}

%%%%%%%%%%%%%%%%%%%%%%%%%%%%%%%%%%%%%%%%%%%%%%%%%%%%%%%%%%%%

%%%%%%%%%%%%%%%%%%%%%%%%%%%%%%%%%%%%%%%%%%%%%%%%%%%%%%%%%%%%
\subsection{Positive elliptic fibrations}
%%%%%%%%%%%%%%%%%%%%%%%%%%%%%%%%%%%%%%%%%%%%%%%%%%%%%%%%%%%%

For holomorphic vector bundles and principal elliptic
bundles, the subject of positivity is well established.
These notions are related geometrically, as we shall
explain now.

Let $X$ be a projective manifold, and $L$ a holomorphic
line bundle on $X$. Positivity of $L$ manifests itself
in the geometric properties of the total space $\Tot(L^*)$
of the non-zero vectors in the dual bundle, which is
a principal $\C^*$-fibration over $X$. When $L$ is positive,
the projection $\Tot(L^*)\stackrel \pi \arrow X$ has no multisections.
Indeed, any such multisection, that is, a subvariety
$X'\subset \Tot(L^*)$ which is surjective and finite over 
$X$, gives a holomorphic section $\gamma_{X'}$ of $(L^*)^{\otimes N}$,
where $N$ is cardinality of $\pi^{-1}(x)\cap X'$ for 
general $x\in X$. The section $\gamma_{X'}$
is defined as a functional on $(L)^{\otimes N}$,
mapping 
\[ v_1 \otimes v_2 \otimes ... \otimes 
   v_N\in L^N\restrict x
\]
to 
\[ \prod\limits_{z_i\in \pi^{-1}(x)\cap X'} v_i(z_i).
\]
When $L$ is positive, $(L^*)^N$ has no sections, hence
$\Tot(L^*)\stackrel \pi \arrow X$ has no multisections.

For example, when $X$ is $\C P^n$ and $L$ is $\calo(1)$,
$\Tot(L^*)$ is identified with $\C^{n+1} \backslash 0$.
This space is quasi-affine and cannot contain compact
subvarieties. 

This effect is even more striking if we take
a quotient $M:=\Tot(L^*)/\langle q \rangle$
of $\Tot(L^*)$ by the $\Z$-action
mapping $v$ to $qv$, where $q\in \C^*$ is
a fixed number, $|q|>1$. This quotient is a 
principal elliptic bundle over $X$, with
the fiber isomorphic to $T= \C^*/\langle q \rangle$.
When $L$ is positive, this manifold is called 
{\bf a regular Vaisman manifold} (see e.g
\cite{_Dragomir_Ornea_}, \cite{_OV:Immersion_}).

 It was shown (see e.g. \cite{_Verbitsky:Sta_Elli_})
positive-dimensional complex subvarieties of
this principal elliptic fibration are $T$-invariant
whenever $L$ is ample. In other words, for any closed subvariety
$Z\subset M$, $Z$ is a preimage of $Z_1\subset X$,
for some complex subvariety of $X$.

If, in addition, $\dim_\C X>1$, a similar result is true
for all stable holomorphic bundles on $M$, and even all
stable reflexive coherent sheaves on $M$.
Namely, given a stable reflexive coherent sheaf
$F$, there exists a line bundle ${\cal L}$ on $M$ such that
$F\cong \pi^* F_1 \otimes {\cal L}$, where $F_1$ is a stable
coherent sheaf on $X$ (\cite{_Verbitsky:Sta_Elli_}).

This result has a number of far-reaching consequences,
even in the most trivial example when $X= \C P^n$,
$L= \calo(1)$ and $M= \C^{n+1} \backslash 0/\langle q\rangle$
is a classical Hopf manifold. However, the argument of
\cite{_Verbitsky:Sta_Elli_} is much more general,
and can be applied to other principal elliptic bundles.

\hfill

%%%%%%%%%%%%%%%%%%%%%%%%%%%%%%%%%%%%%%%%%%%%%%%%%%%%%%%%%%%%
\definition
A {\bf principal elliptic fibration}, or a {\bf principal
  elliptic bundle}
is a complex manifold $M$ equipped with a free holomorphic
action of a 1-dimensional compact complex torus $T$.
We consider a principal elliptic fibration
as a $T$-fibration $M\stackrel \pi\arrow M/T$.

\hfill

Examples of these include the Hopf manifold
and regular Vaisman manifolds defined above, but also
such interesting objects as a Calabi-Eckmann manifold,
that is, a product of two odd-dimensional spheres
$S^{2n-1}\times S^{2m-1}$, $m, n >1$, equipped with an
$SO(2n) \times SO(2m)$-invariant complex structure
(see Subsection \ref{_CE_Subsection_}).
In \cite{_Verbitsky:Sta_Elli_}, the positivity of
principal elliptic bundles was defined as follows.

\hfill

%%%%%%%%%%%%%%%%%%%%%%%%%%%%%%%%%%%%%%%%%%%%%%%%%%%%%%%%%%%%
\definition
Let $M\stackrel \pi\arrow X$ be a principal elliptic
fibration, with $X$ a K\"ahler manifold. We call 
$\pi$ {\bf a positive elliptic fibration}, if
for some K\"ahler form $\omega$ on $X$, $\pi^*\omega$
is exact. 

\hfill

Of course, all principal elliptic fibrations
$M\stackrel \pi\arrow X$ with $H^2(M)=0$ and $X$ K\"ahler
are necessarily positive. This includes all Hopf manifolds
(which are diffeomorphic to $S^1 \times S^{2n-1}$) and
the Calabi-Eckmann manifolds $S^{2n-1}\times S^{2m-1}$.

For a manifold $M:=\Tot(L^*)/\langle q \rangle$ considered
above, positivity is equivalent to the positivity of the
corresponding line bunlde $L$.

The positivity of elliptic fibrations has powerful
geometric applications.

\hfill

%%%%%%%%%%%%%%%%%%%%%%%%%%%%%%%%%%%%%%%%%%%%%%%%
\theorem\label{_positivity_elli_Theorem_}
(\cite{_Verbitsky:Sta_Elli_})
Let $X$ be a compact K\"ahler manifold, and 
$M \stackrel \pi \arrow X$ a positive principal
elliptic fibration. Then, for any positive-dimensional,
irreducible subvariety $Z\subset M$, $Z= \pi^{-1}(Z_1)$,
for some complex subvariety $Z_1\subset X$. If, in
addition, $\dim_\C X>1$, then for any stable
reflexive coherent sheaf $F$ on $M$, 
$F= \pi^* F_1 \otimes {\cal L}$, where $F_1$ is
a stable coherent sheaf on $X$ and $L$ a line bundle.

\endproof

%%%%%%%%%%%%%%%%%%%%%%%%%%%%%%%%%%%%%%%%%%%%%%%%%%%%%%%%%%%%
\subsection{Convexity and positivity}
%%%%%%%%%%%%%%%%%%%%%%%%%%%%%%%%%%%%%%%%%%%%%%%%%%%%%%%%%%%%

For the purposes of the present paper, 
a good notion of positivity in fiber bundles
would give the same consequences as in 
\ref{_positivity_elli_Theorem_}. All 
irreducible subvarieties of a positive fiber
bundle $M \stackrel \pi \arrow X$ should lie 
in the fiber or have the form $\pi^{-1}(Z)$;
and all stable bundles must be pullbacks,
up to tensoring with a line bundle.

Further on in this paper, we study principal
toric fibrations, that is, fibrations with a fiber which
is a compact complex torus. They are defined in the
same way as principal elliptic fibrations.

\hfill

%%%%%%%%%%%%%%%%%%%%%%%%%%%%%%%%%%%%%%%%%%%%%%%%
\definition
Let $T$ be a compact complex torus.
A {\bf principal toric fibration}, or a {\bf principal
  toric bundle} with a fiber $T$
is a complex manifold $M$ equipped with a free holomorphic
action of $T$. We consider a principal toric fibration
as a $T$-fibration $M\stackrel \pi\arrow M/T$.

\hfill

Principal toric fibrations appear in different disguises
in all kinds of algebro-geometric research. 
A classification of these objects is found in
\cite{_Hofer:remarks_}.

\hfill

Even for principal toric fibrations, it is
not clear what is a good notion of positivity. 
One could require that $\pi^*\omega$ is 
exact, for some K\"ahler form $\omega$ on $X$.
This is clearly not enough. Indeed, one could
take a positive elliptic fibration, and 
its product with a torus $T_1$, as a toric
fibration. The pullback $\pi^*\omega$
is exact. However, any subvariety of $X\times T_1$
(and there could be many) gives a corresponding
subvariety of $M$. 

Still, this notion has useful geometric consequences.

\hfill

%%%%%%%%%%%%%%%%%%%%%%%%%%%%%%%%%%%%%%%%%%%%%%%%
\definition
Let $M \stackrel\pi \arrow X$ be a 
fibration of complex manifolds.
We say that $\pi$ is {\bf convex} if for some 
K\"ahler form
$\omega$ on $X$, $\pi^*\omega$ is exact.

\hfill

%%%%%%%%%%%%%%%%%%%%%%%%%%%%%%%%%%%%%%%%%%%%%%%%
\definition
Let $M \stackrel\pi \arrow X$ be a
fibration of complex manifolds, and $X'\subset M$ a
closed subvariety. Assume that the restriction
$\pi\restrict {X'}:\; X' \arrow X$ is proper, surjective
and generically finite. Then $X'$ is called
{\bf a multisection} of $\pi$. 

\hfill

The following trivial lemma explains the utility of
convexity in fiber bundles.

\hfill

%%%%%%%%%%%%%%%%%%%%%%%%%%%%%%%%%%%%%%%%%%%%%%%%
\lemma\label{_convex_multise_Lemma_}
Let $M \stackrel\pi \arrow X$ be a morphism 
of complex manifolds, $X$ compact and K\"ahler,
$\omega$ a K\"ahler class. Assume that $\pi^*\omega$
is exact. Then $\pi$ does not admit multisections.

\hfill

{\bf Proof:} Let $X'$ be a multisection, with 
a generic point $x\in X$ having precisely $N$ pre-images
in $X'$ under the map $\pi\restrict {X'}:\; X' \arrow X$.
Then
\[
\int_{X'} \pi^*\omega^{\dim X} = N \int_X \omega^{\dim X} >0.
\]
This is impossible, because $\pi^*\omega$ is exact.
\endproof

\hfill

The same argument is valid in weaker assumptions.

\hfill

%%%%%%%%%%%%%%%%%%%%%%%%%%%%%%%%%%%%%%%%%%%%%%%%
\definition
Let $M$ be a complex manifold, $\dim_\C M =n$,
and $\omega$ a (1,1)-form, which is positive and closed.
We say that $\omega$ is {\bf weakly K\"ahler} if
$\omega^n\neq 0$ somewhere on $M$. Of course,
in this case $\int_M \omega^n>0$ (assuming $\omega^n$
has compact support). 

\hfill

%%%%%%%%%%%%%%%%%%%%%%%%%%%%%%%%%%%%%%%%%%%%%%%%
\definition
Let $M \stackrel\pi \arrow X$ be a 
fibration of complex manifolds.
We say that $\pi$ is {\bf weakly convex} if for some 
weakly K\"ahler form
$\omega$ on $X$, $\pi^*\omega$ is exact.

\hfill

%%%%%%%%%%%%%%%%%%%%%%%%%%%%%%%%%%%%%%%%%%%%%%%%
\remark\label{_weak_conve_multisec_Remark_}
Clearly, \ref{_convex_multise_Lemma_}
is valid in assumption of weak convexity.

\hfill

Let $M\stackrel \pi\arrow X$ be a principal toric bundle with
fiber $T$, and $V$  a complex variety equipped
with a holomorphic action of $T$. Consider the quotient
$M_V:= (M\times V)/T$, with $T$ acting diagonally on
$M$ and $V$. Clearly, $M_V\stackrel{\pi_V}\arrow X$ is a locally trivial 
holomorphic fibration over $X$, with fiber $V$.

\hfill

%%%%%%%%%%%%%%%%%%%%%%%%%%%%%%%%%%%%%%%%%%%%%%%%
\definition
In this situation, $M_V\arrow X$ is called {\bf a fiber
bundle associated with a principal toric bundle $M\arrow X$}.

\hfill

For principal toric bundles, the following notion
of positivity seems to be useful.

\hfill

%%%%%%%%%%%%%%%%%%%%%%%%%%%%%%%%%%%%%%%%%%%%%%%%
\definition
Let $M\stackrel \pi\arrow X$ be a principal toric bundle with
fiber $T$. Assume that for all surjective homomorphisms
$T \arrow T_1$ from $T$ to a positive-dimensional torus
$T_1$, the associated bundle $M_{T_1}$ is (weakly) convex.
Then $M\stackrel \pi \arrow X$ is called (weakly) Q-positive.

\hfill

The main result of this paper is the following theorem

\hfill

%%%%%%%%%%%%%%%%%%%%%%%%%%%%%%%%%%%%%%%%%%%%%%%%
\theorem\label{_Q_pos_main_Theorem_}
Let $M\stackrel \pi \arrow X$ be a principal toric bundle,
and $Z\subset M$ an irreducible complex subvariety.
Assume that $\pi$ is compact and Q-positive. Then $Z$ is 
$T$-invariant\footnote{That is, takes a form
  $\pi^{-1}(Z_1)$ for some $Z_1\subset X$}, or
$Z$ lies in a fiber of $\pi$. 

\hfill

{\bf Proof:} We prove \ref{_Q_pos_main_Theorem_}
in Section \ref{_rela_Douady_Section_}. \endproof

\hfill

%%%%%%%%%%%%%%%%%%%%%%%%%%%%%%%%%%%%%%%%%%%%%%%%
\remark
For $\dim T=1$, this result follows directly from
\ref{_convex_multise_Lemma_}.

\hfill

%%%%%%%%%%%%%%%%%%%%%%%%%%%%%%%%%%%%%%%%%%%%%%%%
\remark
For $\dim T> 1$, we do not have any control over
the category of coherent sheaves on a positive
$T$-fibration. It is not clear whether 
Q-positivity will imply results similar to
\ref{_positivity_elli_Theorem_}, showing 
that all stable bundles are
pullbacks (up to a multiplication by 
a line bundle). This is the reason for ``Q-''
in Q-positivity (``Q'' stands for ``quotient'').

\hfill

%%%%%%%%%%%%%%%%%%%%%%%%%%%%%%%%%%%%%%%%%%%%%%%%
\remark 
For examples of Q-positive toric fibrations, see
\ref{_Samelson_positive_Theorem_}.

\hfill

%%%%%%%%%%%%%%%%%%%%%%%%%%%%%%%%%%%%%%%%%%%%%%%%
\remark\label{_Q_pos_Subvarieties_Remark_}
Let  $M\stackrel \pi \arrow X$ be a principal toric bundle,
and $X_1\arrow X$ a holomorphic bimeromorphic map.
Then $M\times_X X_1 \stackrel{\pi_1} \arrow X_1$
is a principal toric bundle. When $\pi$ is Q-positive,
$\pi_1$ is not necessarily Q-positive; however,
it is weakly Q-positive. Similarly, if $Z\subset X$ is 
a closed subvariety, $Z_1\arrow Z$ its resolution of
singularities, and 
$M_1:= \pi^{-1}(Z) \times_Z Z_1 \stackrel {\pi_1}\arrow Z_1$
the corresponding principal toric bundle, $\pi_1$ is
weakly Q-positive, though not necessarily Q-positive.

\hfill

Recall that a
compact complex torus $T$ is called {\bf simple} 
if for any holomorphic homomorphism $T\stackrel\phi\arrow T_1$
to another torus, $\ker \psi$ is finite. A generic complex torus
is simple; a generic abelian variety is also
simple. Actually, non-simple tori are quite rare,
and the codimension of the space of non-simple
tori inside the moduli space of all tori grows 
linearly with the dimension of the torus; 
same is true for abelian manifolds. 

\hfill

%%%%%%%%%%%%%%%%%%%%%%%%%%%%%%%%%%%%%%%%%%%%%%%%
\claim \label{_simple_to_Q_posi_Claim_}
Let $M\stackrel \pi \arrow X$ be a principal toric bundle with
fiber $T$. Assume that $T$ is simple and $M$ is (weakly) convex.
Then $M$ is (weakly)  Q-positive. 

\hfill

{\bf Proof:} Let $\omega$ be a (weakly) K\"ahler form on $X$ such
that $\pi^*\omega=d\theta$ is exact, and $T\stackrel \mu\arrow T_1$  a
surjective homomorphism to a positive-dimensional 
torus $T_1$. Since $T$ is simple, $\mu$ is an isogeny
(a finite quotient map). Therefore, 
$M_{T_1}\stackrel{\pi_1}\arrow X$ is a 
quotient of $M$ over a finite group $\Gamma = \ker \mu$
freely acting on $M$.  Averaging $\theta$ with $\Gamma$,
we obtain a 1-form $\theta_\Gamma$ on $M_{T_1}$
which satisfies $d\theta_\Gamma=\pi_1^* \omega$.
\endproof

%%%%%%%%%%%%%%%%%%%%%%%%%%%%%%%%%%%%%%%%%%%%%%%%%%%%%%%%%%%%

\section{Homogeneous complex manifolds}
\label{_homo_mani_Section_}
%%%%%%%%%%%%%%%%%%%%%%%%%%%%%%%%%%%%%%%%%%%%%%%%%%%%%%%%%%%%

%%%%%%%%%%%%%%%%%%%%%%%%%%%%%%%%%%%%%%%%%%%%%%%%
\subsection{Homogeneous complex manifolds and
  toric fibrations}
\label{_homo_mani_toric_Subsection_}
%%%%%%%%%%%%%%%%%%%%%%%%%%%%%%%%%%%%%%%%%%%%%%%%

This paper appeared as an attempt to understand the
algebraic geometry of compact, homogeneous, simply
connected manifolds. These manifolds were classicied
by H. C. Wang (\cite{_Wang:1954_}) and
much studied during 1950-ies and the 1960-ies.
H. Samelson (\cite{_Samelson_}) discovered
homogeneous complex structures on even-dimensional
semisimple Lie groups. In 1980-es and 1990-ies the 
Samelson manifolds became more and more important
in string physics, differential geometry
and algebraic geometry.

Ever since the physicists
Spindel et al (\cite{_SSTvP_}) and then Joyce 
(\cite{_Joyce_}) constructed homogeneous hypercomplex structures
on some compact Lie groups, generalizing Samelson's work,
these manifolds became a favourite testbed for all kinds of conjectures
and constructions in quaternionic geometry
and sometimes in physics (see e.g. \cite{_OP_}). 
In addition, these manifolds are examples
of two rather unique (and mysteriously
related) geometries, which arose in 
string theory. The Joyce's hypercomplex manifolds
are the only known in dimension $>4$ examples 
of so-called {\em strong HKT geometry} (\cite{_Howe_Papado_}, 
\cite{_Gra_Poon_}). In addition, these manifolds
are among the few compact examples of {\em generalized
K\"ahler manifolds} (known in physics as 
{\em manifolds with (2,2)-supersymmetry});
see \cite{_GHR_}, \cite{_Hitchin:biherm_}.
Combining these two structures, one obtains
an even more mysterious (and rare) geometry,
called {\em generalized hyperk\"ahler geometry},
or (4,4)-supersymmetry  (\cite{_GHR_},
\cite{_Bredthauer_}, \cite{_Moraru_Verbitsky_}),
which also occurs on Joyce's hypercomplex manifolds.

\hfill

From their construction, which is explained in Subsection 
\ref{_CE_Subsection_}, 
it is clear that the Samelson manifolds are 
fibered over a flag space, with fibers being compact
tori. It is a special case of a general 
result due to  Borel and Remmert (\cite{_Borel_Remmert_})
and Tits (\cite{_Tits:homoge_}), and its generalization,
proven in this paper (\ref{_Borel_Remmert_Tits_genera_Theorem_});
see also \cite{_Grauert_Remmert:homogene_}. 

\hfill

%%%%%%%%%%%%%%%%%%%%%%%%%%%%%%%%%%%%%%%%%%%%%%%%
\theorem\label{_homo_are_toric_fibr_Theorem_}
Let $M$ be a compact, homogeneous, simply connected
complex manifold. Then $M$ is a principal
toric fibration over a base $X$, which is a homogeneous,
rational, projective manifold.

\hfill

{\bf Proof:} 
Recall that a complex manifold is called {\bf
  parallelizable} if its holomorphic tangent bundle is trivial.
The following lemma is used.

\hfill

%%%%%%%%%%%%%%%%%%%%%%%%%%%%%%%%%%%%%%%%%%%%%%%%
\lemma\label{_triv_K_parall_Lemma_}
(\cite{_Tits:homoge_}) Let $M$, $\dim_\C M=n$ be a compact homogeneous
complex manifold with trivial canonical bundle.
Then $M$ is parallelizable.

\hfill

{\bf Proof:} The holomorphic  tangent bundle $TM$ is
globally generated. Chose a set $\gamma_1, ..., \gamma_n$
of holomorphic vector fields which are linearly
independent at $m\in M$. Since the canonical
bundle of $M$ is trivial, $\gamma_1, ..., \gamma_n$
are linearly independent everywhere. \endproof

\hfill

Return to the proof of
\ref{_homo_are_toric_fibr_Theorem_}.
Since $TM$ is globally generated, the anticanonical
bundle $K^{-1}M$ is globally generated too, hence
base point free. This gives a holomorphic morphism
$M \stackrel\psi \arrow {\Bbb P}(\Gamma (K^{-1}M))$
from $M$ to the projectivization of the space
of sections of $K^{-1}M$. This map is obviously
equivariant with respect to the action of
the group $G= \Aut(M)$. Therefore, the 
projective manifold $B= \psi(M)$ is homogeneous,
and hence rational. By adjunction formula,
the fibers $F$ of $\psi$ have trivial canonical class;
\ref{_triv_K_parall_Lemma_} implies that they
are parallelizable. 
Using the homotopy exact sequence for a fibration
\[
\pi_2 (B) \arrow \pi_1(F) \arrow \pi_1(M)=0,
\]
we obtain that the fundamental group of $F$ is abelian.
Therefore, $F$ is a complex torus (see the proof
of \ref{_Borel_Remmert_Tits_genera_Theorem_} below
for a detailed argument).  The sheaf $T_\psi M$
of holomorphic vector fields tangent to the fibers
of $\psi$ is a trivial bundle, by the same reasining as used to
prove \ref{_triv_K_parall_Lemma_}.
Let $G_0\subset \Aut(M)$ be a Lie group generated by the
sections of $T_\psi M$. Identifying $G$ with the 
group of automorphisms of fibers of $\psi$,
we obtain that $G$ is a compact torus, and
$\psi:\; M \arrow X$ is a principal toric fibration.
\endproof

%%%%%%%%%%%%%%%%%%%%%%%%%%%%%%%%%%%%%%%%%%%%%%%%%%%%%%%%%%%%
\subsection{A generalization of Borel-Remmert-Tits theorem}
%%%%%%%%%%%%%%%%%%%%%%%%%%%%%%%%%%%%%%%%%%%%%%%%%%%%%%%%%%%%

Borel-Remmert-Tits theorem has a generalization, which 
is quite useful in our work.

\hfill

%%%%%%%%%%%%%%%%%%%%%%%%%%%%%%%%%%%%%%%%%%%%%%%%%%%%%%%%%%%%
\theorem\label{_Borel_Remmert_Tits_genera_Theorem_}
Let $M$ be a compact complex manifold equipped with a 
transitive holomorphic action of a Lie group $G$.
Assume that $M$ admits a $G$-invariant decomposition
$M=M_0\times (S^1)^n$, where $\pi_1(M_0)=0$.
Let $\psi:\; M \arrow B\subset {\Bbb P}(H^0(K^{-1}))$
be the $G$-equivarant fibration constructed in the proof of
 \ref{_homo_are_toric_fibr_Theorem_}. Then  $\psi$ 
is a principal toric fibraton.

\hfill

{\bf Proof:} As we have shown, the fibers $F$ of $\pi$ are compact,
homogeneous manifolds which are holomorphically parallelizable,
that is, have trivial tangent bundle. To prove
\ref{_Borel_Remmert_Tits_genera_Theorem_} we have to show that
these fibers are complex tori. It is well known 
(see \cite{_Wolf:parallelism_}) that such manifolds are
always of form $F= {\cal G}/\Gamma$, where ${\cal G}$ is
a complex Lie group, and $\Gamma$ a cocompact lattice. 
To prove that $F$ is a torus, it suffices to show that $\pi_1(F)$ is 
abelian. Indeed, $\Gamma$ is a subgroup of $\pi_1(M)$,
and (being a cocompact lattice), $\Gamma$ is Zariski dense
in ${\cal G}$, hence ${\cal G}$ is also abelian. We 
reduced \ref{_Borel_Remmert_Tits_genera_Theorem_} to the
following lemma.

\hfill

%%%%%%%%%%%%%%%%%%%%%%%%%%%%%%%%%%%%%%%%%%%%%%%%
\lemma\label{_pi_1_fibers_Lemma_}
In assumptions of \ref{_Borel_Remmert_Tits_genera_Theorem_},
we have a decomposition $F=F_0\times (S^1)^n$, where
$\pi_1(F_0) = \pi_2(B)/j(\pi_2(M))$, where 
$j:\; \pi_2(M)\arrow \pi_2(M)$ is the map induced by the
projection.

\hfill

{\bf Proof:} Since $G$ preserves the decomposition
$M= M_0\times (S^1)^n$,  $G$ contains a torus $T^n=(S^1)^n$ acting trivially
on $M_0$ and tautologically on $(S^1)^n$. Consider the natural
action of $T^n\subset G$ on $B$. If any of its orbits has 
positive dimension, it is isomorphic to a quotient
$T= T^n/T_0$. This gives a free holomorphic action
of a torus $T$ on $B$. This is impossible, because
$B$ is rational, and any holomorphic vector field
on $B$ vanishes somewhere. We have shown that $T^n$
acts trivially on $B$, and the map $\psi:\; M \arrow B$
factors through the standard projection
$M = M_0\times (S^1)^n \stackrel \pi \arrow M_0$.
Therefore, $F= F_0 \times T^n$, with $F_0$ being a fiber
of the standard fibration $M_0 \arrow B$. To see that the fundamental
group of $F_0$ is abelian, we use an exact sequence of
homotopy
\[
\pi_2 (M_0)\stackrel j \arrow \pi_2 (B) \arrow \pi_1(F_0) \arrow \pi_1(M_0)=0.
\]
associated to the fibration $M_0\arrow B$. This gives an isomorphism
$\pi_1(F_0) = \pi_2(B)/j(\pi_2(M))$. \ref{_pi_1_fibers_Lemma_}
is proven. We proved \ref{_Borel_Remmert_Tits_genera_Theorem_}.
\endproof

\hfill

%%%%%%%%%%%%%%%%%%%%%%%%%%%%%%%%%%%%%%%%%%%%%%%%
\remark
This theorem has a generalization, also due to J. Tits
(\cite{_Tits:homoge_}). Let $M$ be a compact homogeneous
manifold, with a Lie group $G$ acting on $M$ transitively. 
Then $M$ is fibered over a projective,
rational, homogeneous space with a fiber which
is homogeneous and paralellizable. This fibration
is called {\bf the Tits fibration}. Tits also 
proved that it is $G$-invariant, and universal.

\hfill

Quite often, a compact, homogeneous complex manifold
would have $H^2(M,\Q)=0$. The Calabi-Eckmann manifolds and
the Hopf manifolds we mentioned previously clearly satisfy
$H^2(M)=0$. Let $G$ be a compact Lie group with
$H^1(G, \Q)=0$. The cohomology algebra of $G$ is a
Hopf algebra; by Hopf theorem (\cite{_Hopf:groups_}),
$H^*(G,\Q)$ is a free exterior algebra,
with generators of odd degree. In particular,
$H^2(G,\Q)=0$.

If $G$ is even-dimensional, $G$ admits a left-invariant
complex structure; if $G$ is odd-dimensional,
$G\times S^1$ admits a left-invariant complex structure.
These manifolds are called {\em Samelson spaces}
(see Subsection \ref{_Samelson_spa_Subsection_}).
In both of these cases, $H^2(M, \Q)=0$, and the
corresponding principal toric bundle 
(\ref{_Borel_Remmert_Tits_genera_Theorem_}) 
is convex. If the torus $T$ is, in addition,
simple, $M$ is Q-positive.
Further on in this section, we explain that this
happens for generic complex structures on
$G$ and $G\times S^1$.

%%%%%%%%%%%%%%%%%%%%%%%%%%%%%%%%%%%%%%%%%%%%%%%%%%%%%%%%%%%%
\subsection{The Calabi-Eckmann construction}
\label{_CE_Subsection_}
%%%%%%%%%%%%%%%%%%%%%%%%%%%%%%%%%%%%%%%%%%%%%%%%%%%%%%%%%%%%

The complex structures on the Lie groups can be 
constructed explicitly, in such a way that the
principal toric action becomes apparent.
To see this we first review the construction of 
Calabi-Eckmann manifolds. 

Consider the principal $\C^*\times \C^*$-fibration
\[
  \tilde M = \C^n\backslash 0 \times \C^m\backslash 0 \arrow 
   X = \C P^{n-1} \times \C P^{m-1}.
\]
Given a subgroup $G\subset \C^*\times \C^*$,
we obtain an associated fiber bundle $M_G = \tilde M/G\arrow X$,
with a fiber $\C^*\times \C^*/G$. We parametrize
$\C^*\times \C^*$ by the set of pairs 
$e^{t_1} \times e^{t_2}$, $t_1, t_2 \in \C$. 
Given a number $\alpha \in \C \backslash \R$,
we observe that the subgroup 
\[ G := \{ e^t\times e^{\alpha t} \subset \C^*\times \C^*,
\ \ t\in \C\}\subset \C^*\times \C^*
\]
is co-compact, and the quotient $\C^*\times \C^*/G$ is an
elliptic curve $T$, isomorphic to the quotient
\[
\C / \langle 2\pi\1 n+ 2\pi\alpha \1m\rangle, \ \ n, m \in \Z
\]
Therefore, the fibration $M_G \arrow X$ is an elliptic
fibration, with fiber $T$; it is easy to check that
$M_G$ is actually diffeomorphic to 
$S^{2n-1}\times S^{2m-1}$. We have constructed
the Calabi-Eckmann manifold.

A similar construction will work whenever we have
a principal $(\C^*)^{2n}$-bundle $\tilde M \arrow X$. 
We may identify $(\C^*)^{2n}$ with a quotient of 
$W=\C^{2n}$ by a subgroup 
\[
\Gamma =  \left\{ \sum_{i=1}^{2n} 2\pi\1 k_i\ \ | \ \  k_i \in \Z\right\}
\]
isomophic to  $\Z^{2n}$. Choosing an $n$-dimensional 
complex subspace $V\subset W$ not intersecting 
$\R \Gamma$, we find that the quotient 
$T=W/\langle V+\Gamma\rangle$ is a compact torus
(this is actually a standard construction of 
compact tori, used in their classification).

Now consider $V$ as a subgroup of $(\C^*)^{2n}=\C^{2n}/\Gamma$.
Clearly, the quotient $\tilde M/V$ is a principal toric
bundle over $X$, with fiber $T$. 

%%%%%%%%%%%%%%%%%%%%%%%%%%%%%%%%%%%%%%%%%%%%%%%%%%%%%%%%%%%%
\subsection{Samelson spaces}
\label{_Samelson_spa_Subsection_} 
%%%%%%%%%%%%%%%%%%%%%%%%%%%%%%%%%%%%%%%%%%%%%%%%%%%%%%%%%%%%

The left-invariant complex structures on Lie groups
were constructed by H. Samelson in 1954 (\cite{_Samelson_}).

\hfill

%%%%%%%%%%%%%%%%%%%%%%%%%%%%%%%%%%%%%%%%%%%%%%%%
\definition
Let $G$ be an even-dimensional, compact, semisimple Lie group,
and $I$ a left-invariant complex structure on $G$. Then
$(G,I)$ is called {\bf a Samelson space}. It is a
homogeneous complex manifold.

\hfill

Let $G$ be an even-dimensional
compact Lie group; it is obviously
semisimple. Denote by $G_\C$ its complexification,
and let $B\subset G_\C$ be its Borel subgroup,
and $F= G_\C/B$ the corresponding flag space,
which is obviously compact. The group $G$ can be
considered as a subgroup in $G_\C$; the natural
projection $G\arrow F$ is surjective, and has
the maximal torus of $G$ as its fiber. This is the
toric fibration we require; it remains only
to produce a complex structure on $G$ which
is compatible with this fibration.

Denote by $U$ the
unipotent radical of $B$. The space $\tilde M:=G_\C/U$
is fibered over $F$ with a fiber $T_\C=(\C^*)^k$,
with $k= \dim_\R G - 2 \dim_\C F$ being even.
We arrive at the situation described at the
end of Subsection \ref{_CE_Subsection_}: a $(\C^*)^{2n}$-fibration
over a compact base. Taking a quotient of $\tilde M$
with a subgroup $V\subset T_\C=$ isomorphic to $\C^n$, as indicated in
Subsection \ref{_CE_Subsection_}, we obtain a principal
toric fibration $M= \tilde M/V$. Embedding $G$ 
into $\tilde M= G_\C/U$, and projecting to $M$,
we identify $G$ with $M$, and this identification
is by construction $G$-invariant.

Clearly, this principal fibration is the one
postulated by Borel-Remmert-Tits theorem
(\ref{_homo_are_toric_fibr_Theorem_}).

Notice that we have a complete freedom of choosing
$V\subset T_\C$. Therefore the left-invariant complex structure on 
$G$ can be chosen in such a way that any given torus
of appropriate dimension becomes a fiber of the
Tits fibration.

There is another construction of Samelson manifolds,
which is even more elementary and explicit. Take the Lie algebra $\g$
of a compact, even-dimensional Lie group $G$, and let
$\g_\C$ be its complexification. Take a positive root system 
$e_1, ..., e_n \in \g_\C$, and let $f_1, ..., f_n$ the corresponding
negative roots, and ${\goth h}_\C\subset \g_\C$ be the
corresponding Cartan subalgebra, so that
\[
\g_\C = \langle e_1, ..., e_n,f_1, ..., f_n\rangle \oplus {\goth h}_\C.
\]
Let ${\goth h}$ be a Lie algebra of a maximal torus of $G$
such that ${\goth h}_\C$ is its complexification.
The space $\g$ is decomposed as
\[
\g = {\goth h} \oplus \goth F,
\]
where $\goth F$ is the space of real vectors in
$\langle e_1, ..., e_n,f_1, ..., f_n\rangle$.
Since $\dim G$ is even, the real vector space
${\goth h}$ is even-dimensional.

Now choose a complex structure $\goth I$ on $\g$ , in such a way
that the decomposition $\g = {\goth h} \oplus \goth F$
is preserved. To define a complex structure $\goth I$ on $\goth F$
it is the same as to define a decomposition
\[ 
  {\goth F}\otimes \C= {\goth F}^{1,0}\oplus {\goth
  F}^{0,1},
\] 
with $\overline{{\goth F}^{1,0}}={\goth F}^{0,1}$,
and $\goth I$ acting as $\1$ on ${\goth F}^{1,0}$
and as $-\1$ on ${\goth F}^{0,1}$. We define
\[
  {\goth F}^{1,0}=\langle e_1, ..., e_n \rangle, 
  \ \ {\goth F}^{1,0} =\langle f_1, ..., f_n \rangle.
\]
Choose the action of $\goth I$ on ${\goth h}$ 
arbitrarily. 

Using the left action of $G$ on itself,
$\goth I$ becomes a left-invariant almost
complex structure $I$. By Newlander-Nirenberg,
$I$ is integrable if and only if
$[T^{0,1}G, T^{0,1}G] \subset T^{0,1}G$.
This is equivalent to 
\[ [\g^{0,1}, \g^{0,1}]\subset \g^{0,1},
\]
where $\g^{0,1}\subset \g_\C$ is a space of all
vectors where $\goth I$ acts as $-\1$.
However, from the construction of $\goth I$
it is apparent that $\g^{0,1}$ is identified with
a subspace of Cartan algebra plus
$\langle f_1, ..., f_n\rangle$,
and such a space is closed under commutators.

It is not difficult to check that
this construction is in fact identical
to that we gave earlier. The same argument
was used by Samelson to construct an
invariant complex structure on $G/T$,
where $G$ is a compact Lie group and
$T\subset G$ a compact torus of even 
codimension in $G$.

We conclude this section with the following theorem,
which follows from the constructions we have given.

\hfill

%%%%%%%%%%%%%%%%%%%%%%%%%%%%%%%%%%%%%%%%%%%%%%%%
\theorem\label{_Samelson_positive_Theorem_}
Let $G$ be a compact, even-dimensional Lie group,
and $I$ a left-invariant complex structure on $G$.
Consider the corresponding Tits
fibration $(G,I) \stackrel \pi \arrow F$,
with fibers identified with a compact complex torus $T$. 
Assume that $T$  is simple\footnote{This happens if $I$ is generic,
as shown above.}. Then the fibration $\pi$ is Q-positive.

\hfill

{\bf Proof:} Let \[ G=G_0\times S_1^{k},\] where
$G_0$ is a semi\-simple group with $H^1(G_0, \Q)=0$.
Such a decomposition follows from the classification
of compact Lie groups. It is easy to see that $G$ is a toric
fibration over the flag space $F$ associated with $G_0$
(\cite{_Samelson_}). However, $H^2(G_0, \Q)=0$ as we have
mentioned above. The Tits
map $G\stackrel \pi\arrow F$ is a composition of the projection
$G\arrow G_0$ and the natural map $G_0\arrow F$.
Therefore, the pull-back homomorphism 
\[ \pi^*:\; H^2(F, \R) \arrow H^2 (G, \R)\]
vanishes. We obtain that $\pi^*\omega$ is exact for
any K\"ahler form $\omega$ on $F$. This implies that
$(G,I)$ is convex. Using \ref{_simple_to_Q_posi_Claim_},
we obtain that $(G,I)$ is Q-positive. \endproof

\hfill

%%%%%%%%%%%%%%%%%%%%%%%%%%%%%%%%%%%%%%%%%%%%%%%%
\remark 
The same argument works to show that
for any torus $T\subset G$ of even codimension
in a compact Lie group $G$, a generic left-invariant
complex structure on $G/T$ is Q-positive.
Such spaces were constructed and classified
by Samelson in \cite{_Samelson_}.

%%%%%%%%%%%%%%%%%%%%%%%%%%%%%%%%%%%%%%%%%%%%%%%%%%%%%%%%%%%%

\section{Fiber bundles associated with principal fibrations}

%%%%%%%%%%%%%%%%%%%%%%%%%%%%%%%%%%%%%%%%%%%%%%%%%%%%%%%%%%%%

%%%%%%%%%%%%%%%%%%%%%%%%%%%%%%%%%%%%%%%%%%%%%%%%%%%%%%%%%%%%
\subsection{Varieties of Fujiki class C}
%%%%%%%%%%%%%%%%%%%%%%%%%%%%%%%%%%%%%%%%%%%%%%%%%%%%%%%%%%%%

%%%%%%%%%%%%%%%%%%%%%%%%%%%%%%%%%%%%%%%%%%%%%%%%%%%%%%%%%%%%
\definition
Let $X$, $Y$ be complex varieties, and $C\subset X\times Y$
a closed subvariety, such that the projections
$C\arrow X$, $C\arrow Y$ are proper and generically
one-to-one. Such $C$ is called {\bf a bimeromorphic
morphism}, or {\bf bimeromorphic modification}
between $X$ and $Y$. 

\hfill

Clearly, a composition of bimeromorphic maps is well defined,
and is again bimeromorphic. 

\hfill

%%%%%%%%%%%%%%%%%%%%%%%%%%%%%%%%%%%%%%%%%%%%%%%%%%%%%%%%%%%%
\definition
Let $M$ be a complex variety. We say that $M$ 
{\bf is in Fujiki class C} if $M$ is bimeromorphic
to a K\"ahler manifold. 

\hfill

Varieties of Fujiki class C were defined by A. Fujiki 
in \cite{_Fujiki:1978_}, and explored at great length 
in subsequent papers (\cite{_Fujiki:Douady_78_}, 
\cite{_Fujiki:Douady_82_}). 

The following properties of Fujiki class C will be 
useful for our work.

\hfill

%%%%%%%%%%%%%%%%%%%%%%%%%%%%%%%%%%%%%%%%%%%%%%%%
\proposition\label{_subva_Fujiki_C_Proposition_}
(Fujiki) 
Let $M$ be a variety of Fujiki class C, and $X\subset M$ 
a closed complex subvariety. Then $X$ is also of Fujiki
class C.

\hfill

{\bf Proof:}  \cite{_Fujiki:Douady_78_}, Lemma 4.6
\endproof

\hfill

%%%%%%%%%%%%%%%%%%%%%%%%%%%%%%%%%%%%%%%%%%%%%%%%
\theorem \label{_Douady_Fujiki_C_Theorem_}
Let $M$ be a compact variety of Fujiki class C,
$X\subset M$ a closed subvariety,
and $D$ its Douady deformation space. Then $D$
is compact and belongs to Fujiki class C.

\hfill

{\bf Proof:}  \cite{_Fujiki:Douady_78_}, \cite{_Fujiki:Douady_82_}.
\endproof

\hfill

In \cite{_Demailly_Paun_}, 
 J.-P. Demailly and M. Paun gave a simple characterization
of Fujiki class C.

Recall that a {\bf $p$-current} on a manifold $M$ is a
$p$-form taking values in distributions. 
We consider currents as functionals on smooth
$(\dim_\R M-p)$-forms with compact supports,
called {\bf test-forms}. De Rham differental on
$p$-currents is defined by the formula
\[
\langle d\eta, \psi\rangle := (-1)^{p}\langle \eta, d\psi\rangle,
\]
where $\psi$ is a test-form, and $\eta$ is a $p$-current.
The cohomology of currents are equal to de Rham cohomology
of the manifold. The Hodge decomposition on currents
\[
{\cal C }^p(M) = \bigoplus\limits_{i+j=p} {\cal C }^{i,j}(M)
\]
is defined in the same manner. A (1,1)-current 
$\eta$ is called {\bf positive} if $\langle \eta, \psi\rangle \geq 0$
for any test-form $\psi\in \Lambda^{\dim_\R M-1,\dim_\R M-1}(M)$ 
which is positive, that is, Hodge dual to a pseudo-Hermitian
(1,1)-form with non-negative eigenvalues.

A {\bf closed current} is a current $\eta$ satisfying $d\eta=0$.
A {\bf K\"ahler current} is a closed $(1,1)$-current $\eta$
such that there exists a Hermitian form $\omega$ on $M$
with $\eta-\omega$ is positive. 

The cohomology classes of non-smooth K\"ahler currents
can be non-K\"ahler. For instance, take a divizor $Z\subset M$.
The divizor $Z$ can be considered as a (1,1)-current  by taking
\[
\langle Z, \psi\rangle := \int_Z \psi.
\]
If $Z$ is an exceptional divizor, and $\omega$ a K\"ahler form,
$\lambda\in \R^{\geq 0}$,
the current $\omega+ \lambda Z$ is K\"ahler, but for 
$\lambda$ sufficiently big, its cohomology class is 
clearly non-K\"ahler.

Demailly and Paun proved the following theorem.

\hfill

%%%%%%%%%%%%%%%%%%%%%%%%%%%%%%%%%%%%%%%%%%%%%%%%%%%%%%%%%%%%
\theorem\label{_Kahler_curre_C_Theorem_}
(\cite[Theorem 3.4]{_Demailly_Paun_})
Let $M$ be a compact complex manifold. Then $M$ lies in
Fujiki class C if and only if $M$ supports a K\"ahler current.

\endproof

\hfill

%%%%%%%%%%%%%%%%%%%%%%%%%%%%%%%%%%%%%%%%%%%%%%%%
\remark\label{_Kahler_curre_constr_Remark_}
Such a K\"ahler current is easy to construct. Given
a bimeromorphism $M \sim M_1$, where $M_1$ is K\"ahler,
and using Hironaka desingularization theorem,
we can construct a resolution $M_2$ which dominates
$M$ and $M_1$ and is obtained by sequence of blow-ups with
smooth centers. A blow-up of a K\"ahler manifold
is again K\"ahler, hence $M_2$ is K\"ahler. 
Consider the map $M_2\stackrel \pi \arrow M$; it is 
bimeromorphic and holomorphic. Let $\omega$ be a K\"ahler
form on $M_2$. We consider $\omega$ as a current
$\underline \omega$
on $M$ using the formula
\[
\langle \underline \omega , \psi\rangle := 
\langle \omega, \pi^*\psi\rangle,
\]
This current is obviusly positive and K\"ahler.

\hfill

Further on, we shall need the following proposition.

\hfill

%%%%%%%%%%%%%%%%%%%%%%%%%%%%%%%%%%%%%%%%%%%%%%%%
\proposition\label{_equi_class_G_Proposition_}
Let $M$ be a compact 
complex variety of Fujiki class $C$ and $G$
a compact Lie group acting on $M$ holomorphically.
Then $M$ admits a $G$-equivariant, smooth K\"ahler
resolution of singularities,
with $G$-invariant K\"ahler metric.

\hfill

{\bf Proof:} 
Without restricting the generality, we may assume
that $M$ is smooth. In this case, property C
is equivalent to existence of a K\"ahler current.
Averaging with $G$, we may assume that this
K\"ahler current is $G$-invariant.
The proof of \ref{_Kahler_curre_C_Theorem_}
involves approximating a given K\"ahler current
by K\"ahler currents with controlled logarithmic 
singularities at $Z_k$, and blowing up the corresponding
singular sets (\cite[Theorem 3.2]{_Demailly_Paun_}).
The sets $Z_k$ are obtained explicitly 
as zeroes of multiplier ideal sheaves 
with prescribed weights on $M$
(\cite{_Demailly:Regularization_}). 
These weights can
be chosen $G$-invariant, because (as seen from
the proof) the corresponding weights $\psi_k$
satisfy $|\psi_k - g\psi_k|<C$, for any $g\in G$.
This implies that $Z_k$ can be chosen $G$-invariant. 
Therefore, the K\"ahler blow-up of $M$ can be
achieved by blowing up the $G$-invariant
subvarieties. This implies existence
of a $G$-equivariant K\"ahler resolution
of $M$. Averaging the K\"ahler form with $G$,
we obtain a $G$-invariant K\"ahler metric.
\endproof

%%%%%%%%%%%%%%%%%%%%%%%%%%%%%%%%%%%%%%%%%%%%%%%%%%%%%%%%%%%%
\subsection{Multisections of fiber bundles and positivity}
%%%%%%%%%%%%%%%%%%%%%%%%%%%%%%%%%%%%%%%%%%%%%%%%%%%%%%%%%%%%

Recall that a {\bf multisection} of a fiber bundle
$M_V\stackrel {\pi_V}\arrow X$ is an irreducible
subvariety $X'\hookrightarrow M_V$ which is proper and 
finite over a general point of $X$.

The main result of this section is the following theorem,
which ensures that under appropriate conditions
a fiber bundle associated with a principal toric bundle
does not admit multisections.

\hfill

%%%%%%%%%%%%%%%%%%%%%%%%%%%%%%%%%%%%%%%%%%%%%%%%
\theorem \label{_no_multi_in_fiber_Theorem_}
Let $M\stackrel \pi\arrow X$ be a principal toric
bundle, with fiber $T$, and $V$ a compact complex variety 
of Fujiki class C, equipped with a holomorphic $T$-action.
Consider the associated fiber bundle
$M_V\stackrel {\pi_V}\arrow X$, $M_V = M \times V/T$,
and let $X'\hookrightarrow M_V$ be its multisection. 
Assume that $M$ is compact and weakly Q-positive.
Then $X'$ is $T$-invariant.

\hfill

{\bf Proof:} 

{\bf Step 1.} Without restricting the generality, we may assume that
$X'$ intersects the smooth part of $M_V$. Indeed, if $X'$
lies in the singular part $V_{sing}$ of $V$, we may replace
$M_V$ with $M_{V_{sing}}$. If, again, $X'$ lies in the
singular part of $V_{sing}$, we replace $V_{sing}$ with
its singular part and so on. Using induction by $\dim V$,
we obtain that at some step we arrive at the situation
when $X'$ intersects the non-singular part. 

\hfill

{\bf Step 2.} Let $T_0\subset T$ be a stabilizer $St(z)$
of a generic point $z\in X'$. If $T=T_0$,
\ref{_no_multi_in_fiber_Theorem_} is proven.
Replacing $T$ by $T/T_0$, and $V$ by its subset fixed by
$T_0$, we arrive at assumptions of
\ref{_no_multi_in_fiber_Theorem_}; however, now
$St(z)=0$ for generic $z\in X'$. It suffices to prove
\ref{_no_multi_in_fiber_Theorem_} in this assumption.

\hfill

{\bf Step 3.} Let $V'$ be a $T$-equivariant K\"ahler
resolution of singularities of $V$, which exists by 
\ref{_equi_class_G_Proposition_}. A proper preimage 
of $X'$ under the natural map $M_{V'}\arrow M_V$
is a multisection of $M_{V'}\stackrel {\pi_{V'}}\arrow X$.
Hence it suffices to prove \ref{_no_multi_in_fiber_Theorem_} 
assuming that $V$ is K\"ahler and smooth.

\hfill

{\bf Step 4.} For any
holomorphic action of a torus $T$ on a complex manifold $Z$,
denote by $St_0(z)\subset T$ the connected component of
the stabilizer of $z\in Z$, and by $T_0$ the connected
component of the group 
\[ \{t\in T \ \ | \ \ \text{$t$ acts trivially on Z}\}.
\] Then $St_0(z)=T_0$. Indeed, $St_0(z)/T_0$
acts faithfully on $Sym^*(T_zZ)$, and this gives
a holomorphic embedding of $St_0(z)/T_0$ into
an algebraic group, which is affine, and hence Stein.
Since $St_0(z)$ is holomorphic and compact, such an
embedding is necessarily trivial.

In Steps 1-3, we have reduced
\ref{_no_multi_in_fiber_Theorem_}
to the case when $V$ is K\"ahler and 
smooth, and the action of $T$ is free generically. Because
of Step 4, it remains now to prove \ref{_no_multi_in_fiber_Theorem_} 
assuming that $V$ is K\"ahler and smooth, and for any
$v\in V$, the stabilizer $St(v)$ is finite.

Now \ref{_no_multi_in_fiber_Theorem_}
can be deduced from the following lemma.

\hfill

%%%%%%%%%%%%%%%%%%%%%%%%%%%%%%%%%%%%%%%%%%%%%%%%
\lemma\label{_no_multise_if_free_Lemma_}
Let $M\stackrel \pi \arrow X$ be a principal
toric bundle with fiber $T$, and $V$ -- a K\"ahler
manifold equipped with an action of $T$. Assume that
$X$ is compact, $M$ is weakly convex, and for any $v\in V$, the 
stabilizer $St(v)$ is finite. Consider the associated 
fibration $M_V\stackrel {\pi_V}\arrow X$.
Then $\pi_V$ admits no multisections.

\hfill

{\bf Proof:} Let $\omega$ be a weakly K\"ahler form on $X$
such that 
\begin{equation}\label{_d_theta_omega_Equation_}
\pi^* \omega=d\theta 
\end{equation}
is exact. Clearly, to
show that $\pi_V$ admits no multisections,
it suffices to prove that $\pi_V^* \omega$ is exact
(see \ref{_convex_multise_Lemma_}, \ref{_weak_conve_multisec_Remark_}). 

Averaging $\theta$ over $T$, we may assume that
$\theta$ is $T$-invariant. Denote by $\gotht$ the Lie algebra of $T$.
For all $t\in \gotht$, 
\[
d(\theta\cntrct t) = \Lie_t \theta - (d\theta)\cntrct t=0.
\]
Therefore, $\langle \theta, t\rangle = \theta\cntrct t$ is constant.
Denote by $\theta\check\;\in \gotht^*$ the corresponding form on
$\gotht$. 

Since $T$ is a complex Lie group acting holomorphically,
the natural embedding\footnote{This map is an embedding, because
the  stabilizer $St(v)$ is finite everywhere. This is the only
place we use this assumption in the whole argument.}
 $\gotht \stackrel\rho \hookrightarrow TV$ is $\C$-linear.
Therefore, the K\"ahler form $\omega_V$ of $V$ is non-degenerate
on $\rho(\gotht)$. Since $\omega_V$ is closed and $T$-invariant,
for all $t_1, t_2 \in \gotht$, the function
$\omega_V(t_1,t_2)\in C^\infty V$
is constant (the argument is the same as used above to show that
$\theta\cntrct t$ is constant). This defines a symplectic
form $\omega_V\check\;$ on $\gotht$. Let $t_\theta\in \gotht$
be defined by the formula
\[
\langle \theta\check\;, t\rangle = \omega\check\;(t_\theta, t).
\]
Denote by $\theta_V$ the 1-form $\omega(t_\theta, \cdot)$
on $V$. By Cartan's formula, $\theta_V$ is closed. Moreover,
\begin{equation}\label{_theta_V_t_formula_Equation_}
\langle \theta_V, t\rangle = \theta\check\; (t) = \langle \theta, t\rangle
\end{equation}
for all $t\in \gotht$. 

Consider the projection maps 
\begin{center}
$\xymatrix @C+7mm @R+7mm@!0 {
& M\times V \ar[dl]_{\pi_1} \ar[dr]^{\pi_2} & \\
M &{} & V\\
}$
\end{center}
and let $\Psi:= \pi^*_1 \theta - \pi_2^*\theta_V$. Since
$\theta_V$ is closed, \eqref{_d_theta_omega_Equation_}
gives 
\begin{equation}\label{_d_Psi_pullback_omega_Equation_}
d\Psi = \pi_1^*\pi^*\omega.
\end{equation}
By construction, $\Psi$ is $T$-invariant; from 
\eqref{_theta_V_t_formula_Equation_}, we obtain that
$\langle \Psi, t\rangle=0$ for all $t\in \goth t$,
with $T$ acting on $M\times V$ diagonally. From 
$T$-invariance and $\langle \Psi, t\rangle=0$
it follows that $\Psi$ is a pullback of a form
$\psi$ on $(M\times V)/T = M_V$. Now, \eqref{_d_Psi_pullback_omega_Equation_}
implies that $\pi_V^*\omega= d\psi$. We proved
that $\pi_V^*\omega$ is exact. This finishes the
proof of \ref{_no_multise_if_free_Lemma_} and 
\ref{_no_multi_in_fiber_Theorem_}. \endproof

%%%%%%%%%%%%%%%%%%%%%%%%%%%%%%%%%%%%%%%%%%%%%%%%%%%%%%%%%%%%

\section{Relative Douady spaces}
\label{_rela_Douady_Section_}

%%%%%%%%%%%%%%%%%%%%%%%%%%%%%%%%%%%%%%%%%%%%%%%%%%%%%%%%%%%%

Using \ref{_no_multi_in_fiber_Theorem_},
we can prove \ref{_Q_pos_main_Theorem_} as follows.

Let $M\stackrel \pi \arrow X$ be a Q-positive toric
fibration, and $Z\subset M$ an irreducible complex subvariety,
which does not lie in a fiber of $\pi$.
Denote by $X_1$ the resolution of singularities
of $\pi(Z)$, and let $M_1:= X_1 \times_{\pi(Z)}\pi^{-1}(\pi(Z))\stackrel {\pi_1}\arrow X_1$ be
the corresponding principal toric fibration, and $Z_1$ be
the proper preimage of $Z$ in $M_1$. 
Clearly, $\pi_1$ is a weakly Q-positive fibration
(see \ref{_Q_pos_Subvarieties_Remark_}).

We have reduced \ref{_Q_pos_main_Theorem_} to the following statement.

\hfill

%%%%%%%%%%%%%%%%%%%%%%%%%%%%%%%%%%%%%%%%%%%%%%%%
\proposition\label{_weakly_Q_pos_main_Proposition_}
Let $M\stackrel \pi \arrow X$ be a principal toric bundle, $\dim X>0$,
and $Z\subset M$ an irreducible complex subvariety.
Assume that $\pi$ is weakly Q-positive, and $\pi(Z)=X$.
Then $Z=M$.

\hfill

{\bf Proof:} 
Let $X$ be a complex variety, $Y\subset X$ a compact
subvariety. In this situation, a space
of deformations of $Y$ inside $X$ is defined,
called {\bf the Douady space} (\cite{_Douady_}, 
\cite{_Fujiki:Douady_78_}).\footnote{There are actually
two different definitions of the moduli of deformations
of subvarieties, the Douady space and the Barlet space,
\cite{_Barlet_}. They are bimeromorphically equivalent.
For our purposes any of these two can serve.}
Denote by $D_Y$ the component of the Douady space
containing $Y$. In \cite{_Fujiki:Douady_82_}, A.
Fujiki proved that if $X$ is a compact manifold
of Fujiki class C, same is true for $D_Y$. 

Let $D$ be a Douady space of subvarieties in $T$,
equpped with a natural $T$-action, and $D_\pi:= M_D$ the
associated fiber bundle. We call $D_\pi$ 
{\bf the relative Douady space of $M\stackrel\pi\arrow  X$.}
It is a classifying space of subvarieties of the fibers of $\pi$.

For $x\in X$, let $Z_x:= \pi^{-1}(x) \cap Z$ be the 
subvariety of a torus $\pi^{-1}(x)$ corresponding to $x$.
This gives a multisection $X' \subset D_\pi$,
where $D_\pi$ is the relative Douady space.
Applying \ref{_no_multi_in_fiber_Theorem_},
we find that any multisection of $D_\pi$
is $T$-invariant. This means that $Z_x\subset \pi^{-1}(x)$
is $T$-invariant. However, the only $T$-invariant subvariety
of a torus $T$ is $T$. Therefore, $Z=M$.
\endproof

\hfill

\hfill

\hfill

{\bf Acknowledgements:} I am grateful to J.-P. Demailly
for valuable advice and consultations, and to the referee
for error-checking and simplification of the proof of 
\ref{_no_multi_in_fiber_Theorem_}.

{\small

\hfill

\noindent {\sc 

\noindent
{\sc  Institute of Theoretical and
Experimental Physics,\\
B. Cheremushkinskaya, 25, Moscow, 117259, Russia }

\noindent
\tt verbit@maths.gla.ac.uk, \ \  verbit@mccme.ru 
}% end of small

\end{document}